\documentclass{article}%
\usepackage{amsmath}%
\usepackage{amsfonts}%
\usepackage{amssymb}%

\newtheorem{theorem}{Theorem}

\def\cS{{\mathcal S}}
\def\tr{\mathop{\rm tr}\,}
\def\rank {\mathop{\rm rank}\,}
\def\diag{{\rm diag}}

\def\qed{\hfill\vbox{\hrule width 6 pt
\hbox{\vrule height 6 pt width 6 pt}}}

\begin{document}
\openup .42\jot


\title{\Large Positivity of Partitioned Hermitian Matrices with Unitarily Invariant Norms}

\author{Chi-Kwong Li${}^a$\thanks{E-mails: ckli@math.wm.edu},\;\;
Fuzhen Zhang${}^b$\thanks{E-mails: zhang@nova.edu}\\
\footnotesize{${}^a$ Department of Mathematics,
College of William and Mary}\\
\footnotesize{Williamsburg, VA 23187, USA}
\\
\footnotesize{${}^b$ Division of Math, Science, and Technology,
Nova Southeastern University}\\
\footnotesize{3301 College Ave, Fort Lauderdale, FL 33314, USA}}
\date{}
\maketitle

\centerline{\em In memory of Robert C. Thompson}
\bigskip
\hrule

\bigskip


\noindent {\bf Abstract}

\medskip
We give a short proof of a recent result of Drury on the positivity of
a $3\times 3$ matrix of the form $(\|R_i^*R_j\|_{\tr})_{1 \le i, j \le 3}$
for any rectangular   complex (or real) matrices $R_1, R_2, R_3$ so that
the multiplication $R_i^*R_j$ is compatible for all $i, j$,
where $\|\cdot\|_{\tr}$ denotes the trace norm.
We then give a complete analysis of the problem when the trace norm is
replaced by other unitarily invariant norms.

\medskip
{\footnotesize

}

\bigskip
\hrule
\bigskip

\date{}
\maketitle
\section{Introduction}

Let $A=(A_{ij})_{1\leq i, j\leq m}$ be a partitioned positive semi-definite matrix.
We are interested in the question when
$\big (\|A_{ij}\|\big )_{1\leq i, j\leq m}$ is also positive semi-definite,
where $\|\cdot \|$ is a unitarily invariant norm; see
\cite{MW71, ZFZActa12} and the references therein.

Let $R_1, R_2, R_3$ be rectangular complex matrices
so that the multiplication
$R_j^*R_k$ is compatible for all $(j,k)$ pairs.
Drury \cite{DruryLAA446} recently showed that
the $3\times 3$ matrix $M=(m_{jk})$
is  positive semi-definite if
$$m_{jk}= \tr(|R_j^*R_k|) = \|R_j^*R_k\|_{\tr},$$
where for a (rectangular) matrix $X$,
$X^*$ denotes its conjugate transpose,  $|X|=(X^*X)^{1/2}$ denotes its polar factor,
and  $\|X\|_{\tr}= \tr|X|$ denotes its trace norm,
i.e. the sum of all singular values $s_j(X)$ of $X$.
We will denote by $M_n$ the set of
$n\times n$ complex matrices.

Drury's theorem is a generalization of a result of Marcus and Watkins (see \cite[p.238]{MW71})
asserting that
if $A=(a_{ij})\in M_p$ is positive semi-definite, $1\leq p\leq 3$,   then
so is $(|a_{ij}|)$ (by taking absolute value of each entry). It is known
that for $p\geq 4$, $(|a_{ij}|)$  need not be
positive semi-definite in general.
In this note, we present a short proof of Drury's result.
We then give a complete analysis of the problem when the trace norm is
replaced by other unitarily invariant norms.

\section{Results and proofs}

\begin{theorem} \label{thm1}
Suppose that $A = (A_{ij})_{1 \le i,j \le 3}$ is a partitioned positive semi-definite matrix,
where each $A_{ii}$ is a square matrix,
$i=1, 2, 3$.
Then $(\tr|A_{ij}|)\in M_3$
is also positive semi-definite.
\end{theorem}

\it
\noindent
Proof. \rm We may assume that all $A_{ij}$ are of size $n\times n$ by adding zero rows
(at the bottom) and zero columns (on the right).
Let $A_{12}$ have polar decomposition $P_{12}U$. We may replace
$A$ by $(I_n \oplus U \oplus I_n)A(I_n \oplus U^* \oplus I_n)$ and assume that
$A_{11}, A_{22}, A_{33}, A_{12}=A_{21}$ are positive semi-definite.
Next, in the modified matrix, assume that $A_{23}$ has polar decomposition
$P_{23}V$.
We may replace
$A$ by $(I_n \oplus I_n \oplus V)A(I_n \oplus I_n \oplus V^*)$ and assume that
$A_{11}, A_{22}, A_{33}, A_{12}=A_{21}, A_{23} = A_{32}$ are all positive semi-definite.
Now,  in the resulting matrix, suppose that  $A_{13}$ has polar decomposition
$P_{13}W$. Let $W = X^*DX$,  $D = \diag(\mu_1, \dots, \mu_n)$
for some unitary $X$ and complex units $\mu_1, \dots, \mu_n$.
Then for $\tilde X = X \oplus X \oplus X$, we have
$$\tilde X A \tilde X^* = (\tilde A_{ij}) = \begin{pmatrix}
\tilde P_{11} & \tilde P_{12} & \tilde P_{13} D \cr
\tilde P_{12} & \tilde P_{22} & \tilde P_{23} \cr
D^* \tilde P_{13} & \tilde P_{23} & \tilde P_{33}\cr\end{pmatrix},$$
where $\tilde P_{ij} = X P_{ij} X^*$.
If we remove all the off-diagonal entries of $\tilde P_{ij}$ to  get a diagonal
matrix $Q_{ij}$, then
$$\tr Q_{ij} = \tr \tilde P_{ij} = \tr |\tilde A_{ij}| = \tr |A_{ij}|, \quad
1 \le i, j \le 3.$$

We {claim} that the matrix $(Q_{ij})_{1 \le i, j\leq 3}$ is positive semi-definite.
It will then follow that $(\tr |A_{ij}|) =
(\tr Q_{ij})$ is positive semi-definite; see, e.g. \cite{MarK69, ZFZActa12}.

To prove our {claim},  note that if we take the $(r,r)$ entries of $\tilde A_{ij}$
to form a matrix $\tilde Q_r \in M_3$, then $\tilde Q_r$ is a principal submatrix
of $(\tilde A_{ij})$ and  is positive semi-definite for each $r = 1, \dots, n$.
By the result of $3\times 3$ matrices, we can change the $(1,3), (3,1)$ entries to their
absolute  values to get a positive semi-definite $Q_r \in M_3$.
Because $(Q_{ij})_{1 \le i, j \le 3}$ is permutationally similar to
$Q_1 \oplus \cdots \oplus Q_n$ and thus
is positive semi-definite, we get the desired
conclusion. \qed

The conclusion of Theorem \ref{thm1} may not hold
if we replace the trace norm by other norms on matrices (that can be defined on
$A_{ij}$ for all $1 \le i, j \le n$). For example, suppose
$\|\cdot\|_p$ is the Schatten $p$-norm defined by $\|X\|_p = (\tr|X|^p)^{1/p}$
for $p \in [1,\infty]$, where $\|X\|_{\infty}$ is the operator norm of $X$
and $\|X\|_1$ is the trace norm. Take  the $4\times 4$ positive semi-definite matrix
$A = \left ( {I_2 \atop I_2} {I_2 \atop I_2} \right )$ and partition it into
$(A_{ij})_{1 \le i, j \le 3}$ with $A_{11}=A_{22}=(1)$ and $A_{33}=I_2$. Then
$$
(\|A_{ij}\|_p) = \begin{pmatrix} 1 & 0 & 1 \cr 0 & 1 & 1 \cr 1 & 1 & 2^{1/p}\cr
\end{pmatrix}$$
 is not positive semi-definite if $p > 1$.
So, it is interesting to determine the types of norms $\|\cdot\|$
on matrices such that Theorem \ref{thm1} is valid.
In the following, we
give a complete answer of the problem for
unitarily invariant norms, i.e. norms $\|\cdot\|$ on matrices such that
$\|UAV\| = \|A\|$ for any matrix $A$ and any unitary matrices $U$ and $V$ of
appropriate sizes. One may see \cite{Li} and its references
for some general background of unitarily invariant norms.

Suppose that
$A = (A_{ij})_{1 \le i, j \le m}$ is positive semi-definite.
Through block permutation, we may assume that the diagonal blocks
$A_{11}, \dots, A_{mm}$ have sizes in ascending order.
Let $m'\times m'$ be the size of $A_{mm}$.
Suppose $\|\cdot\|$ is a unitarily invariant norm on $M_{m'}$.
Extend the definition of $\|\cdot\|$ to other blocks by setting
$\|A_{ij}\| = \|\tilde A_{ij}\|$,
where $\tilde A_{ij} \in M_{m'}$ is obtained by adding zero rows
and columns to $A_{ij}$. We can then consider
$(\|A_{ij}\|)$ for any $A = (A_{ij})$.

If $m = 2$, then $(\|A_{ij}\|)$ is positive semi-definite; see \cite{LiMathias}.
If $m \ge 4$,
then we can choose a positive semi-definite matrix $B=(b_{ij}) \in M_m$
such that the matrix $B_0$ obtained by taking the absolute values of the entries
of $B$ is not positive semi-definite; see the example
of Thompson in \cite{MW71}. Let
 $A = (A_{ij})$ such that $A_{ij}$ has $(1,1)$ entry
equal to $b_{ij}$ and all other entries equal to zero.
Then $(\|A_{ij}\|) = \gamma B_0$ is not positive semi-definite, where $\gamma = \|E_{11}\|$
and $\{E_{11}, E_{12}, \dots, E_{nn}\}$ is the standard basis for $M_n$.
For $m = 3$, we have the following.

\medskip
\begin{theorem} \label{thm2}  Consider the set $M(n_1,n_2,n)$ of block matrices in the form
$A=(A_{ij})_{1 \le i, j \le 3}$, where  $A_{11} \in M_{n_1}$, $A_{22} \in M_{n_2}$,
and $A_{33} \in M_{n}$ with $n_1 \le n_2 \le  n$.
Let $\|\cdot\|$ be a unitarily invariant norm on $M_{n}$ and $k = \min\{n_1+n_2,n\}$.
The following conditions are equivalent.
\begin{itemize}
\item[{\rm (a)}]  The matrix $(\|A_{ij}\|)$ is positive semi-definite
whenever $(A_{ij})\in M(n_1,n_2,n)$ is positive semi-definite.
\item[{\rm (b)}]
$\|E_{11} + \cdots + E_{kk}\| = \|E_{11}\|+ \cdots + \|E_{kk}\| = k\|E_{11}\|.$
\end{itemize}
\end{theorem}

\noindent\it Proof. \rm
We may normalize $\|\cdot\|$ so that $\|E_{11}\| = 1$.
By the result in  \cite{Mathias} (see also \cite{Li}),
there is a compact set $\cS$ of real vectors
$v = (v_1, \dots, v_n)$ with $v_1 \ge \cdots \ge v_n \ge 0$ such that
for every matrix $B\in M_n$,
$$\|B\| = \max\{ \|B\|_v: v \in \cS\}
\quad \hbox{ with } \quad
\|B\|_v = \sum_{j=1}^{n} v_j s_j(B),$$
where $s_1(B) \ge \cdots \ge s_n(B)$ are the
singular values of $B$. Because $\|E_{11}\| = 1$, we see that
$$1 = \max\{v_1: (v_1, \dots, v_n) \in \cS\}.$$

Suppose (b) holds. Then
$\cS$ contains a vector $\hat v$ whose
first $k$ components  are equal to 1 so that
$$\|E_{11} + \cdots + E_{kk}\| = \|E_{11} + \cdots + E_{kk}\|_{\hat v}  = k.$$
Hence, for any matrix $B\in M_n$ with rank not larger than $k$, we have
$$\sum_{j=1}^k s_j(B) = \|B\|_{\hat v}
\le \|B\| \le \sum_{j=1}^k \|s_j(B)E_{jj}\| = \sum_{j=1}^k s_j(B).$$
That is, if  $\rank (B)\leq k$, then $\|B\|_{\tr}=\|B\|$.
For $B$ with rank larger than $k$,
$$\sum_{j=1}^k s_j(B) \le \|B\|_{\hat v} \le \|B\|.$$
Let $A = (A_{ij})_{1 \le i, j \le 3} \in M(n_1,n_2,n)$ be positive semi-definite.
If $k = n$, then  $A_{ij}$ has rank at most $k$ and
$\|A_{ij}\| = \|A_{ij}\|_{\tr}$ for all $1 \le i, j \le 3$.
By Theorem \ref{thm1}, $(\|A_{ij}\|)$ is positive semi-definite.
Suppose $n_1 + n_2 \le k < n$. Then  $\rank\!\left ( {A_{13}\atop A_{23}} \right )\leq k$,
and there is a unitary $V\in M_n$ such that the last $n-k$ columns of
$\left ( {A_{13}\atop A_{23}} \right )V$
are zero.
We may replace $A$ by $\tilde A = U^*AU$ without changing
$(\|A_{ij}\|)$, where $U = I_{n_1+n_2} \oplus V$,
and assume that
the last $n-k$ columns of
$\left ( {A_{13} \atop A_{23}} \right )V$
are zero.
Suppose $\hat A = (\hat A_{ij})$ is obtained by deleting the last $n-k$
rows and columns of $A$.
Then $\rank (\hat A_{ij})\leq k$ for every $\hat A_{ij}$.
Because $\hat A_{33} \in M_k$ is a principal
submatrix of $A_{33}$, from the above discussion on $B$, we have
$$\tr|\hat{A}_{33}| \le \sum_{j=1}^k s_j(A_{33})\le  \|A_{33}\|.$$
For  $(i,j) \ne (3,3)$, we have
$$\tr|\hat A_{ij}| = \|\hat A_{ij}\| = \|A_{ij}\|.$$
Thus,
$$(\|A_{ij}\|) =  (\tr|\hat A_{ij}|) + \diag(0,0,\varepsilon),$$
where $\varepsilon = \|A_{33}\| - \tr |\hat A_{33}| \ge 0$.
By Theorem \ref{thm1}, $(\tr|\hat A_{ij}|)$ is positive semi-definite.
It follows that  $(\|A_{ij}\|) $ is positive semi-definite.

\medskip
Suppose (b) does not hold.
Let $s$ be the largest integer such that
$\|E_{11} + \cdots + E_{ss}\| = s\|E_{11}\|$. Then $s < k$.
Choose $\tilde n_1 \le n_1, \tilde n_2 \le n_2$ such that
$\tilde n_1 + \tilde n_2 = s+1 =\tilde n \le n$. Thus, $\tilde n_1, \tilde n_2 \leq s$.
By the choice of $s$, $\|I_{\tilde n_1}\|=\tilde n_1$, $\|I_{\tilde n_2}\|=\tilde n_2$,
and $s\leq \|I_{\tilde n}\|<s+1$. Let $\|I_{\tilde n}\|=s+\delta $ for some
$\delta \in [0, 1)$. Construct the block matrix
$$\tilde A = (\tilde A_{ij})_{1 \le i, j \le 3}
\in M_{\tilde n_1+\tilde n_2 +\tilde n}$$
in which
$$\tilde A_{11} = I_{\tilde n_1}\;\;
\tilde A_{22} = I_{\tilde n_2},\;\;
\tilde A_{33} = I_{\tilde n},\;\;
\tilde A_{21}^t = \tilde A_{12} = O_{\tilde n_1 \times \tilde n_2},$$
$$\tilde A_{31}^t = \tilde A_{13} = \begin{pmatrix}
I_{\tilde n_1} & O_{\tilde n_1 \times \tilde n_2}  \cr \end{pmatrix}, \quad
\tilde A_{32}^t = \tilde A_{23} = \begin{pmatrix}
O_{\tilde n_2 \times  \tilde n_1} & I_{\tilde n_2} \cr \end{pmatrix}.$$
Then $\tilde{A}$ is positive semi-definite and
$$\hat A =  (\|A_{ij}\|)
= \begin{pmatrix} \tilde n_1 & 0 & \tilde n_1 \cr
0 & \tilde n_2 & \tilde n_2 \cr \tilde n_1 & \tilde n_2 & s+\delta\cr
\end{pmatrix}.$$
Observe that
$$\det(\hat A) =\tilde n_1 \tilde n_2(\tilde n_1+\tilde n_2-1+\delta)
-\tilde n_1\tilde n_2^2 - \tilde n_1^2\tilde n_2
= \tilde n_1\tilde n_2(\delta -1) < 0.$$
Hence, $\hat A$ is not positive semi-definite.
If needed,  we can add zero rows and columns to $\tilde A$ to get a matrix $(A_{ij})$ with
$A_{11} \in M_{n_1}, A_{22}\in M_{n_2}, A_{33}\in M_n$.
The matrix $(A_{ij})$ is positive semi-definite while  $\hat{A}$ remains the same.
\qed

\bigskip
\noindent
{\bf \large Acknowledgment}

\smallskip
Li was an affiliate member of the Institute for Quantum Computing, University of
Waterloo, an honorary professor of the Shanghai University, and University of Hong Kong.
The research of Li was supported by the USA NSF and HK RCG.
This project was done while Li and Zhang were participating the 2014
Summer International Program at the Shanghai University.
Li and Zhang would like to thank the support of the Shanghai University.

\end{document}